\documentclass[12pt]{article}
\newtheorem{df}{Definition}[section]
\newtheorem{thm}{Theorem}[section]

\newtheorem{prop}{Proposition}[section]
\newtheorem{rem}{Remark}[section]
\usepackage{amsfonts}
\usepackage{amsmath,amssymb}
\title{A new cohomology theory associated to deformations of Lie algebra morphisms.}
\author {Ya\"el Fr\'egier}
\begin {document}
\date{}
\maketitle

\begin{abstract}
We introduce a new cohomology theory related to deformations of Lie algebra
morphisms. This
notion involves simultaneous deformations
of two Lie algebras and a homomorphism between them.

\end{abstract}

\section { Introduction.}
In his foundatory work \cite{GER} Gerstenhaber developed a theory of deformation of associative and Lie algebras. His theory
links cohomologies of these algebras and a cup-product giving
"obstructions" to deformations. Nijenhuis and
Richardson noticed strong similarities between Gerstenhaber's theory
and the deformations of complex analytic structures on compact
manifolds (\cite{N-R1}). They axiomatized the theory of deformations via the
introduction of graded Lie algebras(\cite{N-R2}). The next step was to
try to find more examples of structures entering under
the scope of those ideas. One such example was given by the theory of
deformations of homomorphisms
(\cite{N-R3}).

 The purpose of this paper is to study another equivalence relation that the
one used in \cite{N-R3} and to introduce a new type of cohomology.
It enables us to deform simultaneously algebras and homorphisms. The
article is organised as follows:

In Section 2 we recall the concept of deformation for Lie algebras and
the results obtained in \cite{N-R3} for the case of homorphisms.
In section 3 we introduce a bundle wich enables a more natural notion of
equivalence.
In section 4 we explore the nature of the cohomology theory that
should be associated with simultaneous deformations. In section 5 we define the complex and give an explicit
formula for the coboundary operator. We check fundamental properties,
and also that, in the particular case when the algebra structures
are fixed, one recovers the classical notions of \cite{N-R3}. Our
approach differs from that of \cite{arn} in which the target Lie algebra
is fixed.

Let us point out that the deformation equation obtained in
this paper cannot be reduced to the form of Maurer-Cartan
equation. This equation is cubic and therefore cannot be
expressed in terms of graded Lie algebras. This aspect is developed in
section 6.

\section { Theoretical background.}
We recall the definition of cohomology of Lie algebras and their
homomorphisms in Section \ref{co}. We then explain how deformations are related to
cohomology. We follow \cite{GER} and \cite{LN} in Section \ref{ln}
for the case of Lie algebras and \cite{N-R1} in Section \ref{nr} for the case of morphisms.

\subsection {Cohomolgy theory of Lie algebras and their homomorphisms}\label{co}
Let $U$ and $V$ be to vector spaces.
Denote $\bigwedge^p(U,V) $ the space of 
p-linear skewsymetric maps from $U$ to $V$ \label{cm}. The direct sum of
these spaces $\bigwedge(U,V)=\bigoplus
_{n \in \mathbb{N}}\bigwedge^p(U,V) $ is naturaly $\mathbb{N}$-graded. 

Consider two structures of Lie algebras $\rho$ on $U$ and  $\theta$ on
$V$. They are not allowed to be deformed. Let $\Phi :(U, \rho)
\longrightarrow (V, \theta)$ be a Lie morphism then $V$ is a
$(U,\rho)$-module via $\Phi$: $u.v := \theta (\Phi (u),v)$. The space
$\bigwedge(U,V)$ is then a complex (Hoschild- Serre cohomology of the
$(U,\rho)$-module $V$ )  whose coboundary operator associated to the triple
($\rho,\theta,\Phi$) is given by

\begin{equation}
\begin {array}{l}
\delta ^p \psi(x_{1},\dots ,x_{p+1}) = \sum_{1\leq  s\leq 
p+1}(-1)^{s}\theta (\Phi(x_{s}) ,А\psi(x_{1},\dots 
,\hat{x_{s}},\dots 
,x_{p+1})) \\ [8pt]
                                  +\sum_{1\leq  s< t 
                                  \leq p+1}(-1)^{s+t-1}А\psi
 (\rho(x_{s},x_{t}),ААx_{1},\dots 
,\hat{x_{s}},\dots ,\hat{x_{t}},\dots ,x_{p+1}). \\
 \end {array}\label{bord}
\end{equation}

 Define $B^p=im\delta^{p-1}$ and $Z^p=ker\delta^{p}$.
One can prove the fundamental property:
$$\delta^{p+1}\circ\delta^p=0 $$
hence $B^p\subset Z^p$.
\begin{rem}{\rm
Please note that $\delta^p,B^p, Z^p$ depend on
$\rho,\theta$ and $\Phi$ but we omit it to simplify the notations.}
\end{rem}

The general case of (\ref{bord}) defines the cohomology associated to the
morphism $\Phi$. Cohomology of the Lie algebra $\rho$ is defined setting
$\rho =\theta$ and $\Phi =identity$ in equation (\ref{bord}).

\subsection {Deformations of algebraic structures : the case of
  Lie algebras.}\label{ln}

Let us recall this classical theory, well described in
\cite{LN}. This will show how one can deduce cohomology formulas from
a geometrical argument which will be usefull while seeking our
formulas. A point $\rho \in \bigwedge^2(U,U)$ is a Lie algebra if it satisfies the
Jacobi identity i.e for all $a,b,c \in U$
\begin{equation}\sum_{\circlearrowright} \rho(\rho(a,b),c)=0 \label{j}
\end{equation}
where $\circlearrowright$ means sum over all cyclic permutations of ${a,b,c}$.

The solutions to the Jacobi identity form an algebraic variety
$\mathcal{L}_{U}$ in
$\bigwedge^2(U,U) $, since introducing coordinates
(structure constants), equation (\ref{j}) becomes a set of quadratic polynomials. Points of this algebraic variety are precisely the Lie algebra
structures on $U$. Given a point $\rho$ in $\mathcal{L}_{U}$, one calls a deformation of $\rho$ a curve
$C(t)$ in $\mathcal{L}_{U}$ such that $C(0)=\rho$. If $C(t)$ is analytic, it can be expanded in series:
$$C(t)=\rho+\sum_{i=1}^{\infty}\rho_i t^i$$
where each $\rho_i \in \bigwedge^2(U,U)$.

 But $C(t)$, obeys Jacoby for
all t,
$$\sum_{\circlearrowright} C(t)(C(t)(a,b),c)=0. $$
Expanding it and identifying the first order terms one gets:
\begin{equation}
\sum_{\circlearrowright} \rho(\rho_1(a,b),c)+\rho_1(\rho(a,b),c)=0
\label{cocalg}
\end{equation}
The left hand side of this equation is noted $\delta^2
\rho_1(a,b,c)$. One can check that $\delta^2
\rho_1 \in \bigwedge^3(U,U)$. Moreover $\delta^2
:\bigwedge^2(U,U)\longrightarrow \bigwedge^3(U,U)$ is linear. It is
the coboundary operator.

One can then conclude that the tangent space to $\mathcal{L}_{U}$ in
$\rho$ is included in the set $Z^2$ of cocycles (i.e solutions $\phi$ of $\delta^2 \phi
=0$). It can be shown that these spaces are in fact the same.

 $GL(U) $Бе has a natural action on  $\mathcal{L}_{U}$
 given by if $g \in GL(U)$ and $\rho \in \mathcal{L}_{U}$ $$(\mathcal{A}_{U}(g)Бе\cdot \rho )(a,b) = g^{-1}\circ \rho 
 (g(a),g(b)) \ \forall a,b \in U \label{aa}.$$

Let $g(t)$ be a one parameter analytic curve in $GL(U) 
$ whose value for $t=0$ is the identity.

Write down the action of this curve on  $\rho$:
\begin{eqnarray*}
       (\mathcal{A}(g(t))Бе\cdot 
\rho) (a,b)& = &(\mathbb{I}-tg_{1}+o(t))\circ\\
             \rho \ 
\Big(a+&tg_{1}(a)&+o(t),b+tg_{1}(b)+o(t)\Big)\\
              & = & \rho (a,b)+t\Big( -g_{1}Бе\big(\rho(a,b)\big) + \rho(
              g_{1}(a),b)+\rho(a, g_{1}(b))Бе\Big) \\
              & &  + o(t)
\end{eqnarray*}
The first order coefficient
defines a tangent vector at $\rho$ to the  $GL(U) $-orbit.
Let us call it $\delta^1 g_1$:
 \begin{equation}\delta^1
g_1(a,b)= -g_{1}Бе\rho(a,b) + \rho( g_{1}a,b)+\rho(a, g_{1}b)Бе\label{cobalg}\end{equation}
$\forall a,b,c \in U$
 Clearly, $\delta^1
g_1 \in \bigwedge^2(U,U)$. Moreover $\delta^1
:\bigwedge^1(U,U)\longrightarrow \bigwedge^2(U,U)$ is linear. It is
the coboundary operator. Its image is the set $B^2$ of
coboundaries. We have
just shown that the set of coboundaries contains the tangent space at $\rho$
to the orbit. The converse is also true and can be shown by
exponentiation.
These two operators $\delta^1$ and $\delta^2$ agree with the
definition (\ref{bord}) and can be viewed as the motivating example.

\subsection {Nijenhuis-Richardson theory of deformations of Lie
  algebra homomorphisms.}\label{nr}

We give a short overview of the theory, mainly to show in which sense
our point of view will differ from the classical approach.

In \cite{GER}, Gerstenhaber deduced the cohomology from the
deformation theory. Nijhenhuis and Richardson then axiomatized his
theoty in \cite{N-R1}. In their point of view a deformation theory is
encoded in a $\mathbb{N}$-graded vector space called a graded Lie
algebra together vith a
coboundary operator. In the case of deformations of Lie algebras, it
coincides with the set of cochains $\bigwedge(U,U)$ and the cobondary
operator as defined in Section $\ref{co}$. The graded vector space has a ``super''-bracket defined as
follows in the case of $\bigwedge(U,V)$:

 if $\phi =\omega
\otimes v \in \bigwedge^p(U,V)\hbox{ and } \psi =\pi \otimes w \in
\bigwedge^q(U,V)$ where $\omega, \pi\in (\bigwedge^p(U))^*$ and
$v,w \in V$  :

$$[\![\phi,\psi ]\!] =\omega \wedge \pi \otimes [v,w]$$ cf (\cite{N-R2}).
The Maurer-Cartan equation takes the form:\begin{equation}
\label{MC}
\delta \phi (t)=-\frac{1}{2} [\![ \phi (t),\phi (t) ]\!]
\end{equation}

To satisfy the axioms, a  deformation theory must also have a group of
symmetry (structure group). It is the exponentiation of the ad
representation of the 0th subspace of the given graded space. It is
linked to the cohomology since coboundaries are then the tangent vectors
to the orbit. In \cite{N-R2}, Nijenhuis and Richardson applied this to
the case of Lie algebra morphisms. The graded vector space seemed more
or less obvious and was given by $\bigwedge(U,V)$. A good candidate
for the coboundary operator was given by equation
(\ref{bord}). Accordingly to the axiomatisation they made, the group
of symmetry had to be $G= \{e^{{\mathrm ad\ } a}, a
\in\bigwedge^0(U,V) \}$. Here ${\mathrm ad\ }$ is understood in the sense of
superbracket ${\mathrm ad\ }a:b\longmapsto [\![ a,b ]\!]$.

 But this
equivalence relation does not seem natural to me: in my opinion, two
morphisms should be equivalent if they are conjugated (change of basis
formula).

\section { The morphism bundle.}

Our idea in the following is to take for starting point an other
equivalence relation. We then deduce the whole deformation structure analogously to what has been
recalled in \ref{ln}.
The most natural notion of equivalence for the linear morphisms from
$U$ to $V$ seems to be:
\begin{df}{\rm Consider the $GL(U) 
\times GL(V)$ action:
\begin{eqnarray*}
     \mathcal{A} : (GL(U) 
\times GL(V)) \times L(U,V) & \longrightarrow & L(U,V)  \\ [8pt]
                   (g,h,\Phi ) & \longmapsto & h^{-1}\circ \Phi \circ g.
\end{eqnarray*}
Two linear maps $\Phi $ and  $\Phi '$ are said to be equivalent if they belong to the
same orbit.} 
\end{df}
Let  $\rho $ be a Lie algebra structure on  $U$ and $\theta $ be a Lie algebra structure on $V$.
One defines the set of Lie algebra morphisms associated to the couple $(\rho ,\theta )$:
$$\mathcal{M}_{\rho ,\theta }=\{ \Phi \in L(U,V) \vert 
\Phi (\rho (a,b)) = \theta (\Phi (a),\Phi (b))  \}Бе$$
The space $\mathcal{M}_{\rho ,\theta }$ is not stable under the
action $\mathcal{A}$.
 Hence we are looking for a natural extension of $\mathcal{A}$ and $
 \mathcal{M}_{\rho ,\theta }$ in order to recover this covariance. 
 
 One has to consider the subbundle 
 $\mathcal{M}$ of the trivial vector bundle $\mathcal{L}_{U}Бе\times 
 \mathcal{L}_{V}Бе\times L(U,V)$ defined by:
 $$\mathcal{M}=\{ (\rho ,\theta ,
 \Phi) \in  \mathcal{L}_{U}Бе\times 
 \mathcal{L}_{V}Бе\times L(U,V) \vert  
 \Phi \in \mathcal{M}_{\rho  ,\theta } \} $$
 where $\mathcal{L}_{U}$ is the algebraic variety of Lie algebra structures on
$U$ (see section \ref{ln}). Differently stated, a point of the base of the bundle $\mathcal{L}_{U}Бе\times 
 \mathcal{L}_{V}$ is the simultaneous data of a Lie algebra structure
on $U$ a Lie algebra structure
on $V$. The fibre over such a couple is constituted by all the Lie
morphisms between these two algebras. Note that $\mathcal{M}$ is not a
 vector bundle. 

 One can then extend the action $\mathcal{A}$ (using notations of \ref{aa})
 \begin{eqnarray*}
     \mathcal{A}: (GL(U) 
\times GL(V)) \times \mathcal{M} & \longrightarrow & \mathcal{M}  \\[8pt]
                   (g,h,(\rho ,\theta ,
 \Phi) ) & \longmapsto & (\mathcal{A}_{U}(g)Бе\cdot \rho 
 ,\mathcal{A}_{V}Бе(h)Бе \cdot \theta  ,h^{-1}\circ \Phi \circ g).
\end{eqnarray*}
The morphism bundle $\mathcal{M}$ is clearly stable with respect to
the action $\mathcal{A}$.

\section {Guessing the formulas from geometry.}

A one-parameter deformation of $(\rho ,\theta ,
 \Phi)$ is a curve $\mathcal{C}_{t}Бе=(\rho _{t}Бе,\theta _{t}Бе,
 \Phi_{t}Бе)$ in $\mathcal{M}$ passing through $(\rho ,\theta ,
 \Phi)$ for  $t=0$. Let us assume in the sequel that $\mathcal{C}_{t}$
 is analytic. This assumption is needed to find the formulas and can
 be relaxed in applications (for instance once one has the formulas,
 one can adopt the framework of formal deformations).

 We will reproduce what has been done in section \ref{ln} in order to
 obtain the cohomology we
 want to define.
  
\subsection {Deducing the 2-cocycles.}
One has to write then the constraints 
  $\mathcal{C}_{t}$ should obey in order to be in $\mathcal{M}$, namely:
  \begin{equation}
\Phi_{t}Бе (\rho_{t}Бе (a,b)) = 
\theta_{t}Бе (\Phi_{t}Бе (a),\Phi 
_{t}Бе(b)) \label{me}    
\end{equation} which is the fibre constraint, and
 \begin{equation}(\rho_{t},\theta_{t}) \subset \mathcal{L}_U \times \mathcal{L}_V
\label{jthe}
\end{equation}  
The study of constraints  (\ref{jthe}) has already
been done in section \ref{ln}.
Developping the fiber constraint (\ref{me}) in series:
\begin{eqnarray*}
\Phi 
(\rho(a,b)+t\rho_{1}Бе(a,b)+o(t) )
+ t\Phi_{1}Бе 
(\rho(a,b)+t\rho_{1}Бе(a,b)+o(t) )+ o(t) = \\
\theta (\Phi(a)+ t\Phi_{1}(a)Бе + 
o(t),\Phi(b)+ t\Phi_{1}(b) +o(t) 
)\\
+t\theta_{1}(\Phi(a)+ t\Phi_{1}(a)Бе + 
o(t),\Phi(b)+ t\Phi_{1}(b) + 
o(t))+o(t)
\end{eqnarray*}\label{dm}
and identifying the first order terms in t one has
$$\Phi_{1}(\rho(a,b))+\Phi(\rho_{1}Бе(a,b))=
\theta(\Phi_{1}(a),\Phi(b))+\theta(\Phi(a),
\Phi_{1}(b)) 
+\theta_{1}(\Phi(a),\Phi(b))$$
which can also be written

\begin{equation}
\begin {array}{c}\underbrace{\Phi_{1}(\rho(a,b))-
\theta(\Phi_{1}(a),\Phi(b))-\theta(\Phi(a),
\Phi_{1}(b))}_{\delta \Phi_{1}Бе}Бе\\
  \underbrace{-\theta_{1}(\Phi(a),\Phi(b))+\Phi
(\rho_{1}Бе(a,b))}_{extra\ term} =0.\\
\end {array}\label {cocy}\end{equation}
One should define the cohomology operator $\Delta^1$ in such a way that the preceding
equation together with equation (\ref{cocalg}) gives the equation of a 2-cocycle:

\begin{equation}\begin{array}{l}\Delta^1\left(  \begin{array}{c} 
                            \rho_1\\
                            \theta_1\\
                           \Phi_1\\
                              \end{array}
 \right)\left(  \begin{array}{l} 
                            \alpha,\beta,\gamma\\
                            x,y,z\\
                           a,b\\
                              \end{array}
 \right)  =  \left(  \begin{array}{l} 
                          \sum_{\circlearrowright} \rho(\rho_1(\alpha,\beta),\gamma)+\rho_1(\rho(\alpha,\beta),\gamma) \\
                           \sum_{\circlearrowright} \rho(\rho_1(x,y),z)+\rho_1(\rho(x,y),z)\\
                          \begin {array}{c}\Phi_{1}(\rho(a,b))-
\theta(\Phi_{1}(a),\Phi(b))-\theta(\Phi(a),
\Phi_{1}(b))\\
  -\theta_{1}(\Phi(a),\Phi(b))+\Phi
(\rho_{1}Бе(a,b)) .\\
\end {array}\\
              \end{array} \right)\\
        \end {array} \label {cocd}
\end{equation}
must be zero.
\subsection {Deducing the 2-coboundaries.}
A coboundary has to be a tangent vector to the orbit.
Let $g(t)\times h(t)$ be a one parameter analytic curve in $GL(U) 
\times GL(V)$ whose value for $t=0$ is the identity.

Write down the action of this curve on  $\Phi$:
\begin{eqnarray*}
        \Phi_{t}Бе & = & \mathcal{A}(g(t),h(t))Бе\cdot 
\Phi \\
              & = & (\mathbb{I}-th_{1}+o(t))\Phi 
(\mathbb{I}+tg_{1}+o(t))\\
              & = & \Phi +t (-h_{1}Бе\Phi + \Phi g_{1}Бе) +o(t) 
\end{eqnarray*}
The first order coefficient $ -h_{1}Бе\Phi + \Phi g_{1}Бе$
clearly defines the projection along the fibre of a tangent vector to the orbit.

This, together with equation (\ref{cobalg}) tells us what shape coboundaries must have:
 
\begin{equation}\begin{array}{l}\Delta^0\left(  \begin{array}{c} 
                            g_1\\
                            h_1\\
                            0\\
                              \end{array}
 \right)\left(  \begin{array}{l} 
                            \alpha,\beta\\
                            x,y\\
                           a\\
                              \end{array}
 \right)  =  \left(  \begin{array}{l} 
                         -g_{1}Бе\rho(\alpha,\beta) + \rho( g_{1}\alpha,\beta)+\rho(\alpha, g_{1}\beta)  \\
                          -h_{1}Бе\theta(x,y) + \theta( h_{1}x,y)+\theta(x, h_{1}y) \\
                     -h_{1}Бе\Phi(a) + \Phi g_{1}(a)
              \end{array} \right)\\
        \end {array} \label {cobd}
\end{equation}

\section {The main definitions and properties.}

In this section we give our main definitions generalizing all the
above examples.

\begin {df}{\rm
The set of p-cochains  $\Lambda^p(U,V)$ is defined by:
$$\textstyle
\Lambda^p(U,V)=\bigwedge^{p+1}(U,U)\bigoplus\bigwedge^{p+1}(V,V)\bigoplus\bigwedge^p(U,V)
$$
where the summands in the right hand side were defined in Section \ref{cm}. We change the 0th order term and take for convention
$\bigwedge^{0}(U,V)=0$.}
\end {df}
Let us introduce the notation $\diamond$: if $\lambda \in
\Lambda^p(V,V)$ and $\Phi \in
\L(U,V)$, one defines $\lambda\diamond\Phi \in
\Lambda^p(U,V)$ by:
$$\lambda\diamond\Phi(x_1,\dots,x_p)=\lambda(\Phi(x_1),\dots,\Phi(x_p))$$
where $x_1,\dots,x_p \in U$
\begin {df}{\rm
Let $X^p=(X_1,X_2,X_3)$ be a p-cochain and let
$(\rho,\theta,\Phi) \in \mathcal{M}$. One defines the coboundary
operator $\Delta^p:\Lambda^p(U,V)\longrightarrow \Lambda^{p+1}(U,V)$
associated to $(\rho,\theta,\Phi)$ (one will forget to mention this
triple in the following):

$$\Delta^p\left(  \begin{array}{l} 
                            X_1\\
                            X_2\\
                           X_3\\
                              \end{array}
 \right)=\left(  \begin{array}{l} 
                           \delta^{p+1}( X_1)\\
                           \delta^{p+1}( X_2)\\
                          \delta^p( X_3)+(-1)^p(\Phi\circ X_1-X_2\diamond\Phi)\\
                              \end{array}
 \right).$$ where  $\delta^{p+1}( X_1)$ and $\delta^{p+1}( X_2)$ are as in
section \ref{ln} and $\delta^p( X_3)$ is given by formula of section \ref{nr}}
\end {df}
\begin{rem}
formulas (\ref{cocd}) and (\ref{cobd}) coincide with this definition.
\end{rem}
\begin {thm}{\rm The following identity is true 
$$\Delta^{p+1}\circ\Delta^p=0.$$ and hence $(\Lambda^p(U,V),\Delta^p)$
forms a complex}.
\end {thm}
Proof

let $X$ be a p-cochain.
By definition of $\Delta^p$:
\begin{eqnarray*} \Delta^{p+1}\circ\Delta^p\left(  \begin{array}{l} 
                            X_1\\[6pt]
                            X_2\\[6pt]
                           X_3\\
                              \end{array}
 \right) & = & \Delta^{p+1}\left(  \begin{array}{l} 
                           \delta^{p+1}( X_1)\\[6pt]
                           \delta^{p+1}( X_2)\\[6pt]
                          \delta^p( X_3)+(-1)^p(\Phi\circ X_1-X_2\diamond\Phi)\\
                              \end{array}
 \right) \\[12pt]
         &  = & \left(  \begin{array}{l} 
                            \delta^{p+2}( \delta^{p+1}( X_1))\\[6pt]
                           \delta^{p+2}( \delta^{p+1}( X_2))\\[6pt]
                        \begin{array}{c}   \delta^{p+1}( \delta^p(
                            X_3)+(-1)^p(\Phi\circ
                            X_1-X_2\diamond\Phi))\\
+(-1)^{p+1}(\Phi\circ\delta^p( X_1)-\delta^p(X_2) \diamond \Phi) \\
\end{array} \\
 \end{array}
 \right) \\[12pt]  
      &  = & \left(  \begin{array}{l} 
                           0\\[6pt]
                           0\\[6pt]
                         \begin{array}{l}
                            \delta^{p+1}((-1)^p(\Phi\circ
                            X_1-X_2\diamond\Phi))\\
+(-1)^{p+1}(\Phi\circ\delta^p( X_1)-\delta^p(X_2)\end{array}\diamond\Phi) \\
                              \end{array}
 \right) \\
\end{eqnarray*}
Hence one has to show that the last component of this cochain vanishes.

Choose any $x_1,\dots,x_{p+2}$ from $U$, then one has:
\begin {eqnarray*}
\delta^{p+1}(\Phi\circ X_1) (x_{1},\dots ,x_{p+2}) & = & \sum_{1\leq  s\leq 
p+2}(-1)^{s}\theta (\Phi(x_{s}) ,А\Phi\circ X_1)(x_{1},\dots 
,\widehat{x}_{s},\dots 
,x_{p+2})) \\
                                  +\sum_{1\leq  s< t 
                                  \leq p+2}(-1)^{s+t-1}А &\Phi&\circ X_1 (\rho(x_{s},x_{t}),ААx_{1},\dots 
,\widehat{x}_{s},\dots ,\widehat{x}_{t},\dots ,x_{p+2}) \\
                                  & = & \Phi(  \sum_{1\leq  s\leq 
p+2}(-1)^{s}\rho (x_{s} ,А X_1(x_{1},\dots 
,\widehat{x}_{s},\dots 
,x_{p+2})) \\
                                   +\sum_{1\leq  s< t 
                                  \leq p+2}(-1)^{s+t-1}А & X_1 & (\rho(x_{s},x_{t}),ААx_{1},\dots 
,\widehat{x}_{s},\dots ,\widehat{x}_{t},\dots ,x_{p+2})) \\
                                  & =  & \Phi(\delta ^p (X_1))  (x_{1},\dots ,x_{p+2}).
\end {eqnarray*}
Hence
$$\delta^{p+1}(\Phi\circ X_1)=  \Phi(\delta ^p (X_1)).$$
Furthermore,
\begin {eqnarray*}
\delta^{p+1}( X_2\diamond\Phi)  (x_{1},\dots ,x_{p+2}) & = & \sum_{1\leq  s\leq 
p+2}(-1)^{s}\theta (\Phi(x_{s}) ,А X_2 (\Phi(x_{1}),\dots 
,\widehat{x}_{s},\dots 
,\Phi(x_{p+2}))) \\
                                   +\sum_{1\leq  s< t 
                                  \leq p+2}(-1)^{s+t-1}А & X_2 &  (\Phi(\rho(x_{s},x_{t})),АА\Phi(x_{1}),\dots 
,\widehat{x}_{s},\dots ,\widehat{x}_{t},\dots ,\Phi(x_{p+2})) \\
                               & =   &     \sum_{1\leq  s\leq 
p+2}(-1)^{s}\theta (\Phi(x_{s}) ,А X_2 (\Phi(x_{1}),\dots 
,\widehat{x}_{s},\dots 
,\Phi(x_{p+2}))) \\
                                   +\sum_{1\leq  s< t 
                                  \leq p+2}(-1)^{s+t-1}А & X_2 & (\theta (\Phi(x_{s}),\Phi(x_{t})),АА\Phi(x_{1}),\dots 
,\widehat{x}_{s},\dots ,\widehat{x}_{t},\dots ,\Phi(x_{p+2})) \\
                                 & = & \delta^{p}( X_2)\diamond\Phi  (x_{1},\dots ,x_{p+2}).         
\end {eqnarray*}
Hence,
$$\delta^{p+1}( X_2\diamond\Phi)=\delta^{p}( X_2)\diamond\Phi.$$
Which completes the proof of the theorem.

\section {The deformation equation}

equation (\ref{MC}) is the main character of deformation theory. It is
known as the Maurer-Cartan (or ``deformation'') equation. The
solutions of this equation precisely correspond to the
deformations. We discuss in this section an analog of this equation in our
framework. We get back to the notations of Section \ref{dm}. If
$\mathcal{C}_t=(\Phi_t,\theta_t,\rho_t)$ is an analytic curve in $\mathcal{M}$,
it can be expanded in series:

$\mathcal{C}_t=\sum_{i=1}^{\infty}\mathcal{C}^it^i$ 
where $ \mathcal{C}^i=\left(
\begin {array}{c}
\rho_i\\
\theta_i\\
\Phi_i\\
\end {array}\right)$.

 Denote $\Delta (\mathcal{C}^i)=\left(
\begin {array}{c}       
\Delta (\mathcal{C}^i)_1\\
\Delta (\mathcal{C}^i)_2\\
\Delta (\mathcal{C}^i)_3\\
\end {array}\right)$
and let us introduce the
notation: $\widetilde{\rho}_t=\rho_t-\rho$,
$\widetilde{\theta}_t=\theta_t-\theta$ and $\widetilde{\Phi}_t=\Phi_t-\Phi$.

\begin {prop}
$\mathcal{C}_t=(\Phi_t,\theta_t,\rho_t)$ is a curve in $\mathcal{M}$ (i.e a
deformation) if and only if $\theta_t$ and $\rho_t$ satisfy equation
(\ref{MC}) and $\Phi_t$ satisfies the equation 
\begin{equation}
\begin {array}{ccc}
\Delta (\mathcal{C}_t)_3 & = &  -\widetilde{\Phi}_t(
\widetilde{\rho}_t) +\theta_t (\widetilde{\Phi}_t(a),\widetilde{\Phi}_t(b))\\
                   & + & \widetilde{\theta}_t (\Phi_t
                   (a),\widetilde{\Phi}_t(b)) +\widetilde{\theta}_t
                   (\widetilde{\Phi}_t(a),\Phi_t(b))
\end {array}\label {defeque}
\end{equation}
\end {prop}

Proof

Identifying nth order terms in equation (\ref{me}) one gets
\begin {eqnarray*}
\Delta (\mathcal{C}^n)_3 (a,b) & = & \sum_{i=1}^{n-1} \Big( -\Phi_i
\rho_{n-i}(a,b) +\theta (\Phi_i(a),\Phi_{n-i}(b)) \\
                    & + & \sum_{j=0}^{n-i} 
 \theta _i(\Phi_j(a),\Phi_{n-i-j}(b)) \Big).
\end {eqnarray*}

One can conclude collecting all these equations in a single formula.

\begin{rem}
This equation (\ref{defeque}) is cubic and therefore cannot be reduced
to the Maurer-Cartan equation (\ref{MC}). 
\end{rem}
\section { Conclusions and outlooks}
Let us stress that our point of view can be
 applied to other cases. For example it can be done for
 associative algebra morphisms and leads to a cohomology formula
 almost identical to the one obtained in this paper.
 
The study  of the right hand side of the
Maurer-Cartan equation (\ref{MC}) led Nijenhuis and Richardson \cite{N-R2} to define the
structure of graded Lie algebras recalled in \ref{nr}. Does there exist an algebraic
structure underlying the right-hand-side of deformation equation (\ref{defeque})?

A fruitfull application in physics of deformation theory is the
concept of deformation quantization of a mechanical system. On the
other hand, Lie morphisms can be found as symetries (moment maps) of
such systems. Since it is possible now to deform simultaneously
Lie algebra structures and morphisms, it is natural to try a
deformation quantization of the symetries of dynamical systems. This problem of conservation of
symetries is known in physics under the name of anomalies. It reduces
to solving a Maurer-Cartan equation ("no ghost" theorem). But this ``deformation equation''' doesn't
come from a deformation. Can our deformation
theory help to fill this gap?

\vskip 1cm

\noindent
{\bf Acknowledgments.} I am  grateful to V.Ovsienko for  fruitful
discussions and motivations. I would also like to thank C.Duval for
his constant help and his kindness, C.Roger, F.Pellegrini and J.Merker
for their interest on this work and valuable remarks.

\vskip 1cm


\noindent
Ya\"el FREGIER\\
{\small C.N.R.S., C.P.T.}\\
{\small  Luminy-Case 907}\\
{\small  F-13288 Marseille Cedex 9, France}\\
\&\\
{\small Institut
Girard Desargues, URA CNRS 746}
\\ {\small Universit\'e Claude Bernard - Lyon I}\\
{\small 43
bd. du 11 Novembre 1918}\\
{\small 69622 Villeurbanne Cedex, France}

\end{document}